\documentclass[11pt,a4paper]{article}

\usepackage{amsmath}
\usepackage{amssymb,latexsym}
\usepackage{amsthm}
\usepackage{enumerate}
\usepackage[hang,flushmargin,bottom]{footmisc}
\usepackage[margin=0.35in,format=hang,font=small,labelfont=bf]{caption}
\usepackage[normalem]{ulem}
\usepackage{needspace}
\usepackage{graphicx,tikz}
\usepackage[colorlinks=true,urlcolor=blue,linkcolor=blue,citecolor=blue]{hyperref}
\usepackage{palatino}
\usepackage{eulervm}

\newtheorem{thmO}{Theorem}
\newtheorem{corO}[thmO]{Corollary}

\newcommand{\+}{\hspace{0.07 em}}
\newcommand{\thalf}{\tfrac{1}{2}}
\newcommand{\nequiv}{\not\equiv}
\newcommand{\liminfty}[1][n]{\lim\limits_{#1\rightarrow\infty}}
\newcommand{\floor}[1]{\left\lfloor #1 \right\rfloor}
\newcommand{\ceil}[1]{\left\lceil #1 \right\rceil}
\newcommand{\bbZ}{\mathbb{Z}}
\newcommand{\zr}{\bbZ_r}
\newcommand{\zt}{\bbZ_t}
\newcommand{\ztr}{\bbZ_t^{\,r}}

\newcommand{\mybox}[2][0.9851]{\fbox{\parbox{#1\linewidth}{#2}}}
\newcommand{\thmbox}[1]{\mybox{\vspace{-7pt}#1\vspace{-7pt}}}

\newenvironment{bullets} {\vspace{-9pt}\begin{itemize}\itemsep0pt} {\end{itemize}\vspace{-9pt}}

\title{\textbf{Large butterfly Cayley graphs and digraphs}}

\author{$\phantom{{}^\dagger}$David Bevan${}^\dagger$}

\date{}

\begin{document}
\maketitle

{\let\thefootnote\relax\footnotetext
{${}^\dagger$Department of Computer and Information Sciences, University of Strathclyde, Glasgow, Scotland.}}

{\let\thefootnote\relax}

\begin{abstract}
\noindent
We present
families of
large
undirected and directed
Cayley graphs 
whose construction is related to butterfly networks.
One approach yields,
for every large~$k$ and for values of~$d$ taken from a large interval,
the largest known Cayley
graphs and digraphs
of diameter~$k$ and
degree~$d$.
Another
method
yields, for sufficiently large~$k$ and infinitely many values of~$d$,
Cayley
graphs and digraphs
of diameter~$k$ and
degree~$d$ whose order is exponentially
larger in~$k$
than any previously constructed.
In the directed case, these are within a linear factor in $k$ of the Moore bound.
\end{abstract}

\section{Introduction}
The goal of the \emph{degree--diameter problem} is to determine the largest possible order of
a graph or digraph,
perhaps restricted to some special class,
with given maximum (out)degree and diameter.
For an overview of progress on a wide variety of approaches to this problem, see the
survey by Miller \& {\v{S}}ir{\'a}{\v{n}}~\cite{MS2013}.

Our concern here is with 
large
\emph{Cayley} graphs and digraphs.
Recall that, for a group $G$ and a unit-free generating subset $S$ of $G$, the \emph{Cayley digraph} of $G$ generated by $S$ has vertex set $G$ and a directed edge from $g$ to $gs$ for all $g\in G$ and $s\in S$.
If $S$ is symmetric, i.e. $S=S^{-1}$, then
the corresponding undirected simple graph is the \emph{Cayley graph} of $G$ generated by $S$.
The Cayley (di)graph is thus regular of (out)degree $|S|$ 
and vertex-transitive.

We are interested in graphs and digraphs
of 
degree~$d$ and diameter $k$, 
for arbitrary large $k$ and varying $d$.
If a construction yields graphs of order $n_{d,k}$, we say that it 
has \emph{asymptotic order} $f(d,k)$ if, 
for fixed $k$,
$$
\liminfty[d]
\frac{n_{d,k}}{f(d,k)}
\;=\;
1
.
$$

No graph or digraph can be larger than the corresponding \emph{Moore bound}.
For undirected graphs, this bound is
$\mathrm{M}_{d,k}=1+ \frac{d}{d-2}\big((d-1)^k-1\big)$ if $d>2$.
In the directed case, it is
$\mathrm{DM}_{d,k}=\frac{1}{d-1}\big(d^{k+1}-1\big)$ if $d>1$.
In both cases, the Moore bound has asymptotic order $d^k$.


Previously, for arbitrary degree and diameter,
the largest known directed Cayley graphs were obtained by Vetr\'ik~\cite{Vetrik2012} and Abas \& Vetr\'ik~\cite{AV2017}, whose constructions have asymptotic order
$k\big(\frac{d}{2}\big)^{\!k}$ for even $k$, and $2k\big(\frac{d}{2}\big)^{\!k}$ for odd $k$.
Our construction yields Cayley digraphs whose order is asymptotically $k\+d^{k-1}$.
For fixed diameter $k\geqslant8$, these digraphs are larger than those in~\cite{Vetrik2012} and~\cite{AV2017} for every value of $d$ in a large interval.
We also construct,
for fixed $k$ and infinitely many values
of $d$,
Cayley digraphs whose asymptotic order is
$\frac{d^k}{e^2k}$, 
a factor of $\frac{2^{k-1}}{e^2k^2}$ larger than those of Abas \& Vetr\'ik, and
within a linear factor in $k$ of the Moore bound.

The undirected case is similar. Previously, the largest known Cayley graphs
were obtained by Macbeth, {\v{S}}iagiov{\'a}, {\v{S}}ir{\'a}{\v{n}} \& Vetr\'ik~\cite{MSSV2010}, whose construction has asymptotic order $k\big(\frac{d}{3}\big)^{\!k}$.
For $d-k\nequiv 3\!\pmod{4}$,
we construct Cayley graphs whose order is asymptotically $k\big(\frac{d}{2}\big)^{\!k-1}$.
For sufficiently large diameter~$k$, these graphs are larger than those in~\cite{MSSV2010} for every suitable value of $d$ in a large interval.
We also construct,
for given $k$ and infinitely many values
of $d$,
Cayley graphs whose asymptotic order is
$\frac{1}{e^2k}\big(\frac{d}{2}\big)^{\!k}$, 
a factor of $\frac{1}{e^2k^2}\big(\frac{3}{2}\big)^{\!k}$ larger than those
in~\cite{MSSV2010}.

Our constructions are based on a two-parameter family of groups.
For $t\geqslant2$, let $\zt=\bbZ/t\bbZ$ be the additive group of integers modulo $t$, and for $r\geqslant2$, let $\ztr$ denote the product $\zt\times\ldots\times\zt$, where $\zt$ occurs $r$ times,
considered as an additive group of vectors.
Let $\alpha$ be the automorphism of $\ztr$, defined by $\alpha(v_0,\ldots,v_{r-1})=(v_{r-1},v_0,\ldots,v_{r-2})$,
that cyclically shifts coordinates rightwards by one, and
consider the semidirect product $G=\ztr\rtimes\zr$, of order $r\+t^r$,
with the group operation given by $(u,s)\!\cdot\!(v,s')=(u+\alpha^s(v),s+s')$,
for $u,v\in\ztr$ and $s,s'\in\zr$.
We write elements of $G$ in the form $(v_0,\ldots,v_{r-1};s)$, where each $v_i\in\zt$ and $s\in\zr$. 
Using this notation, the group operation is
\vspace{-3pt}
\begin{multline*}
(u_0,\ldots,u_{r-1};\,s) \!\cdot\! (v_0,\ldots,v_{r-1};\,s') \\
\;\,=\;\,
(u_0+v_{r-s},\, 
\ldots,\, u_{s-1}+v_{r-1},\, u_s+v_0,\,
\ldots,\, u_{r-1}+v_{r-1-s};\, s+s' ) ,
\end{multline*}
arithmetic in the subscripts being performed modulo $r$.
The group $G$ is used to create all our Cayley graphs and digraphs.

The Cayley digraph generated by elements of $G$ of the form $(a,0,\ldots,0;1)$, $a\in\zt$ is isomorphic to the base-$t$ order-$r$ (wrapped) \emph{butterfly network}, $B_t(r)$,
so called because it is composed 
of
$rt^{r-1}$ edge-disjoint
\mbox{$t$-\emph{butterflies}} (copies of the complete bipartite graph $K_{t,t}$); see~\cite[Figure~2]{ABR1990}.
Butterfly networks are closely related to the \emph{de~Bruijn graphs}~\cite{deBruijn1946},
the directed base-$t$ order-$r$ de~Bruijn graph being a coset graph of
$B_t(r)$~\cite[Theorem~4.4]{ABR1990}.

Cayley graphs and digraphs of $G$ were used previously by
Macbeth, {\v{S}}iagiov{\'a}, {\v{S}}ir{\'a}{\v{n}} \& Vetr\'ik~\cite{MSSV2010}
and Vetr\'ik~\cite{Vetrik2012} 
in the constructions mentioned above, though in neither case is
the connection to the
butterfly networks
made explicit.
Each of our results is a consequence of choosing an appropriate set of generators for $G$.
We make use of two distinct constructions.

\Needspace*{7\baselineskip}
\section{The first construction}

We present the directed case first, since it is slightly simpler.

\thmbox{
\begin{thmO}\label{thm1stDir}
For
any $k \geqslant 4$ and $d \geqslant k-1$, there exist Cayley digraphs that have
diameter $k$, outdegree~$d$, and order $(k-1)(d-k+3)^{k-1}$.
\end{thmO}
}

\begin{proof}
Let $r=k-1$ and $t=d-k+3$, and let the underlying group of the Cayley digraph be $G=\ztr\rtimes\zr$.
The order of $G$ is $r\+t^r=(k-1)(d-k+3)^{k-1}.$

As generators for the Cayley digraph 
we use
the $t$ \emph{shift and add} elements 
$(a,0,\ldots,0;1)$, for each $a\in\bbZ_t$,
together with the remaining $r-2$ nonzero \emph{cyclic shift} elements 
$(0,\ldots,0;s)$, for $2\leqslant s\leqslant r-1$.
Thus the digraph has outdegree $t+r-2=d$.

It also has diameter $r+1=k$. Every element
is the product of
$r$ shift and add operations (establishing the vector) and possibly a single cyclic shift (to establish the final shift value if it is nonzero).
On the other hand, if $s\neq0$ then $(1,\ldots,1;s)$ cannot be obtained as a product of fewer than $k$ generators.
\end{proof}

Clearly, the butterfly network $B_t(r)$ is a subdigraph of the Cayley digraph of Theorem~\ref{thm1stDir}.
The additional edges in our construction, corresponding to the cyclic shift elements, consist of $t^r$ vertex-disjoint copies of the complete digraph on $r$ vertices with a directed $r$-cycle removed.

Vetr\'ik~\cite{Vetrik2012} presents, for any $k\geqslant3$ and $d\geqslant4$,
a family of Cayley digraphs
of diameter~$k$, degree~$d$, and order $k\!\floor{\frac{d}{2}}^{\!k}$.
For 
odd diameters,
Abas \& Vetr\'ik~\cite{AV2017} improve this result by a factor of two, constructing Cayley digraphs
of diameter at most~$k$ and degree~$d$ of order $2k\!\floor{\frac{d}{2}}^{\!k}$.
Clearly, for large enough $d$, these digraphs are bigger than those of Theorem~\ref{thm1stDir}.
However,
for any given diameter $k\geqslant8$,
it can be confirmed (using a computer algebra system, or otherwise) 
that the digraphs of Theorem~\ref{thm1stDir} are larger than those of Vetr\'ik and Abas \& Vetr\'ik if
$$
2\+k + 2\ln k 
\;\;<\;\; d \;\;<\;\;
2^{k-1}\big(1-\tfrac{1}{k}\big) -k^2. 
$$
For specific values of the degree, we can do much better.
If we set $d=k^2-3k$, then the digraphs of Theorem~\ref{thm1stDir} have orders
at least $\mathrm{DM}_{d,k}/ek$,
within a linear factor of the Moore bound, and
exceeding those of Abas \& Vetr\'ik by a factor of at least $2^{k-1}/{e k^2}$, which exceeds $1$ for $k\geqslant9$.

For the undirected case, we simply add elements to
the generating set to
make it symmetric. 

\thmbox{
\begin{thmO}\label{thm1stUndir}
For any $k \geqslant 5$ and $d \geqslant k$ such that $d-k\nequiv 3\!\pmod{4}$,
there exist Cayley graphs that have diameter $k$, degree $d$, and order $(k-1)\big(\!\floor{\frac{d-k}{2}}+2\big)^{\!k-1}$.
\end{thmO}
}
\begin{proof}
Let $r=k-1$ and $t=\floor{\frac{d-k}{2}}+2$, and let 
$G=\ztr\rtimes\zr$.
As generators for the Cayley graph of $G$ we use
the $t$ elements $(a,0,\ldots,0;1)$, along with their inverses $(0,\ldots,0,-a;-1)$,
and the remaining $r-3$ nonzero elements $(0,\ldots,0;s)$ for $2\leqslant s\leqslant r-2$.
In addition, if $d-k\equiv 1\!\pmod{4}$, in which case $t$ is even, then the involution $(0,\ldots,0,\frac{t}{2};0)$ is also included as a generator.

Thus the graph has degree $2t+r-3 +(d-k\!\mod2)= d$.
As in the directed case,
it has diameter $r+1=k$.
Every element
is the product of
$k-1$ shift and add operations and possibly a single cyclic shift.
On the other hand, if $s\notin\{-1,0,1\}$ then $(1,\ldots,1;s)$ cannot be obtained as a product of fewer than $k$ generators,
and $G$ has such an element since $r\geqslant4$.
\end{proof}
Macbeth, {\v{S}}iagiov{\'a}, {\v{S}}ir{\'a}{\v{n}} \& Vetr\'ik~\cite{MSSV2010} present, for any $k \geqslant 3$ and $d \geqslant 5$,
a family of Cayley graphs with diameter at most $k$, degree $d$, and
order no greater than 
$k\big(\frac{d+1}{3}\big)^{\!k}$.\footnote{
The graphs in~\cite{MSSV2010} are slightly larger than those of Macbeth, {\v{S}}iagiov{\'a} \& {\v{S}}ir{\'a}{\v{n}}~\cite{MSS2012}, whose order is at most $k\big(\frac{d+1}{3}\big)^{\!k}-k$.}
Their constructions also use the group $G$, with a different generating set.
For large enough $d$, these graphs are bigger than those of Theorem~\ref{thm1stUndir}.
However,
for $k\geqslant27$, 
the graphs of Theorem~\ref{thm1stUndir} are larger than those of Macbeth, {\v{S}}iagiov{\'a}, {\v{S}}ir{\'a}{\v{n}} \& Vetr\'ik
for any $d-k\nequiv 3\!\pmod{4}$ satisfying
$$
 3\+k + 6\ln k 
 \;\;<\;\; d \;\;<\;\;
 2\big(\tfrac{3}{2}\big)^{\!k}\big(1-\tfrac{1}{k}\big)-k^2 . 
$$
For specific values of the degree, we can do much better.
If we set $d=k^2-2k$, then the graphs of Theorem~\ref{thm1stUndir} have orders
exceeding those in~\cite{MSSV2010} by a factor of at least $\frac{2}{e k^2}\big(\tfrac{3}{2}\big)^{\!k}$, which exceeds $1$ for $k\geqslant14$.


\section{The second construction}

In our second construction, we conceive of the 
vectors of length $r$
as being partitioned into $k-1$ \emph{long} blocks, each of length $\ell$, and a single \emph{short} block, of length $m$.

Again, the directed case is presented first, since it is simpler.

\thmbox{
\begin{thmO}\label{thm2ndDirConstruct}
For any $k,\ell,t\geqslant2$ and positive $m<\ell$, there exist Cayley digraphs that have
diameter $k$, outdegree~$t^\ell+(r-1)t^m-1$, and order $r\+t^r$, where $r=(k-1)\ell+m$.
\end{thmO}
}

\begin{proof}
As before, let
$G=\ztr\rtimes\zr$, of order
$r\+t^r$.
As generators for the Cayley digraph, 
we use the $t^\ell$ \emph{long} elements $(a_1,\ldots,a_\ell,0,\ldots,0;\ell)$, $a_i\in\zt$,
together with the additional $(r-1)t^m-1$ nonzero \emph{short} elements $(a_1,\ldots,a_m,0,\ldots,0;s)$, $a_i\in\zt$, $s\neq\ell$.
Thus the digraph has outdegree $t^\ell+(r-1)t^m-1$.
Long elements shift by $\ell$ and modify a long block; short elements shift arbitrarily and modify a short block.

The digraph has diameter $k$. 
Every element 
is the product of
a single short element (establishing $m$ components of the vector and guaranteeing the final shift value)
and
$k-1$ long elements (establishing the remaining $(k-1)\ell=r-m$ components of the vector).
On the other hand, $(1,\ldots,1;0)$ cannot be obtained as a product of fewer than $k$ generators.
\end{proof}
The Cayley digraph of Theorem~\ref{thm2ndDirConstruct}
contains both of the butterfly networks $B_{t^\ell}(r)$ and $B_{t^m}(r)$ as subdigraphs.
Its edges can be partitioned into
$rt^{r-\ell}$ copies of the
\mbox{$t^\ell$-butterfly}, from the long elements,
$r(r-2)t^{r-m}$ copies of the
\mbox{$t^m$-butterfly}, from the short elements that have nonzero shift,
and a collection of directed cycles from the short elements with zero shift.

Given $k$, $\ell$ and $t$, for judicious choice of $m$, these digraphs are larger than those of 
Abas \& Vetr\'ik~\cite{AV2017}.
For example,
if we let $t=2$, then for all $k\geqslant31$ and sufficiently large $\ell$, the order of our digraphs
is greater than that of those in~\cite{AV2017}
if
$$
 \ell - k - \log_2\! \ell + 2
 \;\;<\;\;
 m
 \;\;<\;\;
 \ell - \log_2\! k\ell -\tfrac{2}{k}(\log_2\! k + 2)
 .
$$

If $m$ is chosen optimally, we can do much better than that.

\thmbox{
\begin{corO}\label{cor2ndDir}
  For any $k \geqslant 3$, there are arbitrarily large values of $d$ for which there exist Cayley digraphs that have diameter $k$, outdegree $d$, and order at least $\tfrac{1}{k}\big(\frac{k}{k+2}(d+1)\big)^{\!k}$.
\end{corO}
}

\begin{proof}
  We use the construction of Theorem~\ref{thm2ndDirConstruct}.
  Let $t=2$, and let $\ell$ be any sufficiently large positive integer such that $\log_2k^2\ell\leqslant\tfrac{3}{4}\ell$.
  Let $r=\ceil{k\ell - \log_2 k^2\ell}$, and $m=r-(k-1)\ell$, so $r=(k-1)\ell+m$.
  Note that $0<m<\ell$.

The digraph has diameter $k$ and order $r \+ 2^r$, which (rounding $r$ down) is at least
$$
n_0
\;=\;
\big(k\ell - \log_2 k^2\ell\big)2^{k\ell - \log_2 k^2\ell}
\;=\;
\left(\frac{1}{k} - \frac{\log_2k^2\ell}{k^2\ell}\right)2^{k\ell} .
$$
Its degree is $d=2^\ell+(r-1)2^m-1$, which (substituting for $m$ and rounding $r$ up) is less than
$$
d^+
\;=\;
2^\ell \:+\: \big(k\ell - \log_2 k^2\ell\big) 2^{k\ell - \log_2 k^2\ell + 1 - (k-1)\ell} \:-\: 1
\;=\;
\left(1+\frac{2}{k}-\frac{2\log_2k^2\ell}{k^2\ell}\right)2^\ell \:-\: 1 .
$$

Let $\theta = \frac{\log_2k^2\ell}{k\ell}$.
Note that the condition on $\ell$ implies that $\theta \leqslant \frac{3}{4k}\leqslant\frac{1}{4}$, since $k\geqslant3$.

Now,
$$
k \+ n_0 \left(\frac{k}{k+2}(d^++1)\right)^{\!-k}
\;=\;
\left(1 - \theta\right)
\left(1 + \frac{2\theta}{k+2-2\theta} \right)^{\!\!k}
\;>\;
\left(1 - \theta\right)
\left(1 +  \frac{2k\theta}{k+2-2\theta} \right)
,
$$
which is at least 1 if $k\geqslant2$ and $0\leqslant\theta\leqslant\frac{k-2}{2k-2}$.
Since $k\geqslant3$ and $\theta \leqslant \frac{1}{4}$, the result follows.
\end{proof}
These digraphs have asymptotic order exceeding
$\tfrac{d^k}{e^2k}$,
a factor of $\tfrac{2^{k-1}}{e^2k^2}$ larger than those of Abas \& Vetr\'ik, and within a linear factor in $k$ of the Moore bound.

It is worth briefly explaining the choice of values for $t$ and $r$ in the proof of Corollary~\ref{cor2ndDir}.
Suppose we fix $t$ and $r$ (and hence the order $rt^r$), and also fix the diameter $k$.
What is the optimal choice for $\ell$, that minimises the degree $t^\ell+(r-1)t^{r - (k - 1)\ell}-1$?
Differentiating with respect to $\ell$ and equating to zero yields $\ell=\frac{1}{k}\big(r+\log_t(k - 1) (r - 1)\big)$.
Solving for $r$ then gives
$$
r \;=\;
\frac{1}{\ln t}
\+
W\!\left(\frac{t^{k \ell-1} \ln t}{k-1}\right) \:+\: 1,
$$
where $W$ is the \emph{Lambert W function}, defined implicitly by $W(z)e^{W(z)}=z$. Asymptotically, $W(z) = \ln z-\ln\ln z+o(1)$. Applying this approximation for $W$ then yields $r\approx k\ell - \log_t\! k^2\ell$.
Using this value for $r$ results in a digraph whose order is asymptotically at least
$\tfrac{1}{k}\big(\frac{k}{k+t}(d+1)\big)^{\!k}$. Setting $t=2$ makes this maximal.

The results in the undirected case are similar.
As before, we just add elements to
the generating set to
make it symmetric. 

\thmbox{
\begin{thmO}\label{thm2ndUndirConstruct}
For any $k,\ell,t\geqslant2$ and positive $m<\ell$, there exist Cayley graphs that have
diameter $k$, degree~$2t^\ell + (2 r - 3) t^m - r$, and order $r\+t^r$, where $r=(k-1)\ell+m$.
\end{thmO}
}

\begin{proof}
Let 
$G=\ztr\rtimes\zr$.
As generators for the Cayley graph of $G$ with these parameters, we use:
\begin{bullets}
\item the $t^\ell$ long elements  $(a_1,\ldots,a_\ell,0,\ldots,0;\ell)$, $a_i\in\zt$
\item their $t^\ell$ inverses  $(0,\ldots,0,a_1,\ldots,a_\ell;-\ell)$
\item the $(r-2)(t^m-1)$ short elements  $(a_1,\ldots,a_m,0,\ldots,0;s)$, $a_i\in\zt$ not all zero, $s\notin\{0,\ell\}$
\item their $(r-2)(t^m-1)$ inverses  $(0,\ldots,0,\overbrace{a_1,\ldots,a_m,0,\ldots,0}^s;-s)$
\item the $t^m-1$ nonzero short elements  $(a_1,\ldots,a_m,0,\ldots,0;0)$; this set is symmetric 
\item the $r-3$ short elements  $(0,\ldots,0;s)$, $s\notin\{0,\pm\ell\}$; this set is also symmetric 
\end{bullets}
Thus the graph has degree $2t^\ell + (2 r - 3) t^m - r$.
As in the directed case, it has order $r\+t^r$ and diameter $k$.
\end{proof}

Given $k$, $\ell$ and $t$, for appropriate choice of $m$, these graphs are larger than those of Macbeth, {\v{S}}iagiov{\'a}, {\v{S}}ir{\'a}{\v{n}} \& Vetr\'ik~\cite{MSSV2010}.
For example,
if we let $t=2$, then for all $k\geqslant69$ and sufficiently large $\ell$, the order of our graphs
is greater than that of those in~\cite{MSSV2010}
if
$$
 \ell +k - \log_2\! 3^k\ell + 1
 \;\;<\;\;
 m
 \;\;<\;\;
 \ell - \log_2\! k\ell -\tfrac{3}{k}(\log_2\! k + 2) - 1
 .
$$

If $m$ is chosen optimally, we have the following.

\thmbox{
\begin{corO}\label{cor2ndUndir}
  For any $k \geqslant 3$, there are arbitrarily large values of $d$ for which there exist Cayley graphs that have diameter $k$, degree $d$, and order at least
  $$
  \tfrac{1}{k}\Big(\tfrac{k}{2k+4}\big(d \:+\: k\log_2 \tfrac{d}{2} \:-\: \log_2\log_2 d \:-\: \log_2 8k^2 \big)\!\Big)^{\!k}
  .
  $$
\end{corO}
}

\begin{proof}
We use the construction of Theorem~\ref{thm2ndUndirConstruct}.
As in the proof of Corollary~\ref{cor2ndDir},
  let $t=2$, and let~$\ell$ be any sufficiently large positive integer such that $\log_2k^2\ell\leqslant\tfrac{3}{4}\ell$.
  Let $r=\ceil{k\ell - \log_2 k^2\ell}$, and $m=r-(k-1)\ell$, so $r=(k-1)\ell+m$. 

The graph has diameter $k$ and order $r \+ 2^r$, which is at least
$$
n_0
\;=\;
\big(k\ell - \log_2 k^2\ell\big)2^{k\ell - \log_2 k^2\ell}
\;=\;
\left(\frac{1}{k} - \frac{\log_2k^2\ell}{k^2\ell}\right)2^{k\ell} .
$$
Its degree is $d=2^{\ell+1} + (2 r - 3) 2^m - r$, which (substituting for $m$ and rounding $r$ up in the second term) is less than
$$
  2^{\ell+1} \:+\: \big(2k\ell - 2\log_2 k^2\ell - 1\big) 2^{k\ell - \log_2 k^2\ell+1-(k-1)\ell} \:-\: r \\
  \;=\; 
  \left(2+\frac{4}{k}-\frac{1+4\log_2k^2\ell}{k^2\ell}\right)2^\ell \:-\: r
  .
$$
Thus, $\thalf(d+r)$ is less than
$
q = 
\left(1+\frac{2}{k}-\frac{2\log_2k^2\ell}{k^2\ell}\right)2^\ell ,
$
and by the argument in the proof of Corollary~\ref{cor2ndDir} (with $q=d^++1$), we know that
$
k n_0
> 
\big(\frac{kq}{k+2}\big)^{\!k}
>
\big(\frac{k}{2k+4}(d+r)\big)^{\!k}
.
$

It remains to establish the appropriate lower bound for $r$.

Now, $k n_0<2^{k\ell}$ and $q>\tfrac{d}{2}$, so $2^\ell>\frac{k d}{2k+4}$ and thus
$
\ell
> 
\log_2\tfrac{k d}{2k+4}
= 
\log_2 \tfrac{d}{2}-\log_2\big(1+\tfrac{2}{k}\big)
.
$

Since
$\big(1+\frac{2}{k}\big)^k < e^2 < 2^3$, we have
$\log_2\big(1+\tfrac{2}{k}\big)<\frac{3}{k}$ and thus
$\ell > \log_2 \tfrac{d}{2}-\frac{3}{k}$.

Now, $r \geqslant {k\ell - \log_2 k^2\ell}$, so
$$
r
\;>\;
 k\log_2 \tfrac{d}{2} \:-\: 3 \:-\: \log_2k^2 \:-\: \log_2 \!\big(\log_2 \tfrac{d}{2} -\tfrac{3}{k}\big)
,
$$
which 
is greater than
$
k\log_2 \tfrac{d}{2} - \log_2\log_2 d - \log_2 8k^2 ,
$
as required.
\end{proof}
These graphs have asymptotic order exceeding
$\tfrac{1}{e^2k}\big(\tfrac{d}{2}\big)^{\!k}$,
a factor of $\tfrac{1}{e^2k^2}\big(\tfrac{3}{2}\big)^{\!k}$ larger than those of Macbeth, {\v{S}}iagiov{\'a}, {\v{S}}ir{\'a}{\v{n}} \& Vetr\'ik.

\subsubsection*{Acknowledgements}

The author is very grateful to Grahame Erskine for pointing out errors in earlier drafts, 
and also to an assiduous referee whose detailed comments greatly improved the quality of this paper.

\emph{S.D.G.}

\bibliographystyle{plain}
{\footnotesize\bibliography{../bib/mybib}}

\end{document}